\documentclass{conm-p-l}
\usepackage{amssymb, amsthm, amsmath, amscd}

\newtheorem{theorem}{Theorem}[section]
\newtheorem{lemma}[theorem]{Lemma}

\newtheorem{cor}[theorem]{Corollary}

\theoremstyle{definition}
\newtheorem{example}[theorem]{Example}

\theoremstyle{remark}
\newtheorem{remark}[theorem]{Remark}

\def\red{\mathrm{red}}
\def\cd{\operatorname{cd}}
\def\rad{\operatorname{rad}}
\def\ara{\operatorname{ara}}
\def\Spec{\operatorname{Spec}}
\def\Proj{\operatorname{Proj}}
\def\depth{\operatorname{depth}}

\def\a{\mathfrak a}
\def\b{\mathfrak b}
\def\p{\mathfrak p}
\def\m{\mathfrak m}
\def\et{{\acute{e}t}}

\def\AA{\mathbb A}
\def\CC{\mathbb C}
\def\NN{\mathbb N}
\def\PP{\mathbb P}
\def\QQ{\mathbb Q}
\def\ZZ{\mathbb Z}

\def\calH{\mathcal H}
\def\calO{\mathcal O}

\begin{document}

\title{On the arithmetic rank of certain Segre products}

\author{Anurag K. Singh}
\address{School of Mathematics, Georgia Institute of Technology, Atlanta,
GA~30332.}
\email{singh@math.gatech.edu}
\thanks{The first author is supported in part by the NSF under Grant
DMS~0300600.}

\author{Uli Walther}
\address{Department of Mathematics, Purdue University, West Lafayette, 
IN~47907.}
\email{walther@math.purdue.edu}
\thanks{The second author is supported in part by the NSF under Grant
DMS~0100509}



\begin{abstract}
We compute the arithmetic ranks of the defining ideals of homogeneous
coordinate rings of certain Segre products arising from elliptic curves. The
cohomological dimension of these ideals varies with the characteristic of the
field, though the arithmetic rank does not. We also study the related
set-theoretic Cohen-Macaulay property for these ideals.
\end{abstract}
\maketitle

In \cite{Ly1} Lyubeznik writes: {\em Part of what makes the problem about the
number of defining equations so interesting is that it can be very easily
stated, yet a solution, in those rare cases when it is known, usually is highly
nontrivial and involves a fascinating interplay of Algebra and Geometry.}

In this note we present one of these rare cases where a solution is obtained:
for a smooth elliptic curve $E \subset \PP^2$, we determine the arithmetic rank
of the ideal $\a$ defining the Segre embedding $E \times \PP^1 \subset \PP^5$,
and exhibit a natural generating set for $\a$ up to radical. The ideal $\a$ is
not a set-theoretic complete intersection and, in the case of characteristic
zero, we use reduction modulo $p$ methods to prove moreover that $\a$ is not
set-theoretically Cohen-Macaulay.

\section{The Segre embedding of $E \times \PP^1$}

Let $A$ and $B$ be $\NN$-graded rings over a field $A_0=B_0=K$. The {\em Segre
product}\/ of $A$ and $B$ is the ring
$$
A\# B=\bigoplus_{n\ge 0} A_n\otimes_K B_n
$$ 
which is a subring, in fact a direct summand, of the tensor product
$A\otimes_K B$. The ring $A\# B$ has a natural $\NN$-grading in which
$(A\# B)_n=A_n\otimes_K B_n$. If $U \subset \PP^r$ and $V \subset \PP^s$ are
projective varieties with homogeneous coordinate rings $A$ and $B$
respectively, then the Segre product $A\# B$ is a homogeneous coordinate ring
for the Segre embedding $U \times V \subset \PP^{rs+r+s}$.

Let $E$ be a smooth elliptic curve over a field $K$. Then $E$ can be embedded
in $\PP^2$, and so $E=\Proj \left(K[x_0,x_1,x_2]/(f)\right)$ where $f$ is a
homogeneous cubic polynomial. The Segre product $E \times \PP^1$ has an
embedding in $\PP^5$ with homogeneous coordinate ring
$$
S=\frac{K[x_0,x_1,x_2]}{(f)} \,\,\# \,\, K[y_0,y_1],
$$
i.e., $S$ is the subring of $K[x_0,x_1,x_2,y_0,y_1]/(f)$ generated by the six
monomials $x_iy_j$ for $0 \le i \le 2$ and $0 \le j \le 1$. The relations
amongst these generators arise from the equations
$$
x_iy_j \cdot \ x_ky_l=x_iy_l \cdot \ x_ky_j
$$
defining $\PP^2 \times \PP^1 \subset \PP^5$, and the four equations
$$
y_0^{3-k} y_1^k f(x_0,x_1,x_2)=0 \quad \text{for} \quad 0 \le k \le 3.
$$
The ring $S$ is normal since it is a direct summand of a normal ring. The local
cohomology of $S$ with support at its homogeneous maximal ideal $\m$ may be
computed by the K\"unneth formula, which shows that $H^2_\m(S)$ is isomorphic
to
$$
[H^2_\m(S)]_0=H^1(E,\calO_{E}) \otimes H^0(\PP^1,\calO_{\PP^1})=K.
$$
In particular, $S$ is a normal domain of dimension $3$ which is not
Cohen-Macaulay.

Let $R=K[z_{ij}: 0\le i \le 2, \ 0 \le j \le 1]$ be a polynomial ring. Then $R$
has a $K$-algebra surjection onto $S$ where
$$
z_{ij} \mapsto x_iy_j.
$$
We use $f_{x_i}$ to denote the partial derivative of $f(x_0,x_1,x_2)$ with
respect to $x_i$. The Euler identity implies that
$$
3 f(x_0,x_1,x_2)=\sum_{i=0}^2 x_i f_{x_i}(x_0,x_1,x_2),
$$
and multiplying by $y_0^2 y_1$ and $y_0 y_1^2$ we obtain, respectively, the
equations
$$
3 y_0^2 y_1 f(x_0,x_1,x_2)=\sum_{i=0}^2 z_{i1}f_{x_i}(z_{00},z_{10},z_{20})
$$
and
$$
3 y_0 y_1^2 f(x_0,x_1,x_2)=\sum_{i=0}^2 z_{i0}f_{x_i}(z_{01},z_{11},z_{21}).
$$
Consequently if $K$ is a field of characteristic other than $3$, then the
kernel of the surjection $R \longrightarrow S$ is the ideal $\a$ generated by
the seven polynomials $z_{10}z_{21}-z_{20}z_{11}$, $z_{20}z_{01}-z_{00}z_{21}$,
$z_{00}z_{11}-z_{10}z_{01}$, and
$$
\sum_{i=0}^2 z_{is}f_{x_i}(z_{0t},z_{1t},z_{2t}),
\quad \text{where} \quad 0 \le s,t \le 1.
$$
We prove that the ideal $\a$ has arithmetic rank four, and that the last four
polynomials above generate $\a$ up to radical:

\begin{theorem}\label{thm-ara}
Let $E=\Proj \left(K[x_0,x_1,x_2]/(f)\right)$ be a smooth elliptic curve and
let $\a \subset R=K[z_{ij}: 0\le i \le 2, \ 0 \le j \le 1]$ be the defining
ideal of the homogeneous coordinate ring of the Segre embedding
$E \times \PP^1 \subset \PP^5$. If the characteristic of the field $K$ does
not equal $3$, then the arithmetic rank of the ideal $\a$ is $4$ and
$$
\a=\rad \left( R\sum_{i=0}^2 z_{is}f_{x_i}(z_{0t},z_{1t},z_{2t}) : 
0 \le s,t \le 1 \right).
$$
\end{theorem}

\begin{proof}
Let $(p_0,p_1,p_2,q_0,q_1,q_2) \in \AA^6$ be a root of the four polynomials
above which, we claim, generate $\a$ up to radical. Then
\begin{center}
$f(p_0,p_1,p_2)=0$, \\
$q_0f_{x_0}(p_0,p_1,p_2)+q_1f_{x_1}(p_0,p_1,p_2)+q_2f_{x_2}(p_0,p_1,p_2)=0$, \\
$p_0f_{x_0}(q_0,q_1,q_2)+p_1f_{x_1}(q_0,q_1,q_2)+p_2f_{x_2}(q_0,q_1,q_2)=0$, \\
$f(q_0,q_1,q_2)=0$.
\end{center}
We claim that the size two minors of the matrix
$$
\begin{pmatrix}
p_0 & p_1 & p_2 \\
q_0 & q_1 & q_2
\end{pmatrix}
$$
must be zero. If not, then $P=(p_0,p_1,p_2)$ and $Q=(q_0,q_1,q_2)$ are distinct
points of $\PP^2$. The first and fourth equations imply that the points $P$ and
$Q$ lie on the elliptic curve $E \subset \PP^2$. The second equation implies
that $Q$ lies on the tangent line to $E$ at point $P$, and similarly the third
equation implies that $P$ lies on the tangent line to $E$ at $Q$. But then the
secant line joining $P$ and $Q$ meets the elliptic curve with multiplicity $4$,
which is not possible. Alternatively, consider the group law on $E$ using an
inflection point as the identity element of the group. Then the second and
third equations imply that $2P+Q=0=P+2Q$ in the group law, and so $P$ and $Q$
are the same point of $\PP^2$. This proves the claim, and it follows that
$(p_0,p_1,p_2,q_0,q_1,q_2) \in \AA^6$ is a zero of all polynomials in the ideal
$\a$. By Hilbert's Nullstellensatz, $\a$ is generated by the four polynomials
up to radical once we tensor with $\overline{K}$, the algebraic closure of $K$.
Since $\overline{K}$ is faithfully flat over $K$, the same is true over $K$ as
well and, in particular, $\ara \a \le 4$.

Since the ideal $\a$ has height $3$, the proof of the theorem will be complete
once we show that $\ara \a \neq 3$. In the next section we use \'etale
cohomology to prove, more generally, that the defining ideal of any projective
embedding of $E\times \PP^1$ is not a set-theoretic complete intersection,
Theorem~\ref{etale}. However, in characteristic $0$, we can moreover show that
the local cohomology module $H^4_\a(R)$ is nonzero, from which it follows that
$\ara \a \ge 4$. (In positive characteristic, $H^4_\a(R)=0$ if the elliptic
curve $E$ is supersingular.) We proceed with the characteristic zero case, and
our main tool here is the connection between the local cohomology modules
supported at $\a$, and topological information about the affine variety defined
by $\a$.

If $\ara \a=3$, then each of the $7$ generators of $\a$ has a power belonging
to an ideal $(g_1,g_2,g_3)R \subseteq \a$. This gives us seven equations, each
with finitely many coefficients, and hence we may replace $K$ by the extension
of $\QQ$ obtained be adjoining these finitely many coefficients. A finitely
generated field extension of $\QQ$ can be identified with a subfield of $\CC$,
and so it suffices to prove the desired result in the case $K=\CC$.

Let $U$ be the complement of $X=E\times\PP^1$ in $\PP^5$ and let $V\subset
\AA^6$ be the cone over $U$. The de Rham cohomology of $V$ can be computed from
the \v Cech-de Rham complex \cite[Chapter II]{Bott-Tu} corresponding to any
affine cover of $V$. The de Rham functor on an $n$-dimensional affine smooth
variety, when applied to any module, can only produce cohomology up to degree
$n$ (in our case, $n=5$). It follows that the de Rham cohomology of $V$ will be
zero beyond the sum of $n$ and the index $t$ of the highest nonvanishing local
cohomology module $H^t_\a(R)$.

On the other hand, there is a Leray spectral sequence,
\cite[Thm.~14.18]{Bott-Tu},
$$
H^i_{dR}(U;\calH^j_{dR}(\CC^*;\CC)) \Longrightarrow H^{i+j}_{dR}(V;\CC)
$$
corresponding to the fibration $V\to U$ with fiber $\CC^*$. Since $U$ arises
through the removal of a variety of codimension two from the simply connected
space $\PP^5$, it follows that $U$ is simply connected as well. Hence in the
above spectral sequence the local coefficients are in fact constant
coefficients. The nonzero terms in the $E_2$-page of the spectral sequence have
$j=0,1$. Hence if $k$ denotes the index of the top nonzero de Rham cohomology
group of $V$, then $k-1$ is the index of the top nonzero de Rham cohomology
group of $U$ and $H^{k-1}_{dR}(U;\CC)\cong H^{k}_{dR}(V;\CC)$.

We now claim that $H^8_{dR}(U;\CC)$ in nonzero. Note that
$H^1_{dR}(E;\CC)\cong\CC^2$, so the K\"unneth formula gives
$H^1_{dR}(X;\CC)\cong\CC^2$. Since sheaf and de Rham cohomology agree, using
Alexander duality \cite[V.6.6]{Iversen} and the compactness of $X$, we obtain
$H^{2\cdot 5-1}_X(\PP^5;\CC)\cong \CC^2$. There is a long exact sequence of
sheaf cohomology
$$
\begin{CD}
\cdots @>>> H^i(U;\CC) @>>> H^{i+1}_X(\PP^5;\CC) @>>> H^{i+1}(\PP^5;\CC)
@>>> \cdots,
\end{CD}
$$
which, since $H^9(\PP^5;\CC)=0$, implies that $H^8(U;\CC)\neq 0$. It follows
that $H^9(V;\CC)$ is nonzero as well, and so the local cohomology module
$H^4_\a(R)$ must be nonzero.
\end{proof}

\section{Elliptic curves in positive characteristic}

Let $R$ be a ring of prime characteristic. We say that $R$ is {\em $F$-pure\/}
if the Frobenius homomorphism $F:R \longrightarrow R$ is pure, i.e., if 
$F \otimes 1_M : R \otimes_R M \longrightarrow R \otimes_R M$ is injective for
all $R$-modules $M$. By \cite[Proposition 6.11]{HR}, a local ring $(R,\m,K)$ is
$F$-pure if and only if the map
$$
F \otimes 1_{E_R(K)} : R \otimes_R E_R(K) \longrightarrow R \otimes_R E_R(K) 
$$
is injective where $E_R(K)$ is the injective hull of the residue field $K$.

Let $E$ be a smooth elliptic curve over a field $K$ of characteristic $p>0$.
The Frobenius induces a map
$$
F: H^1(E,\calO_E) \longrightarrow H^1(E,\calO_E) 
$$
on the one-dimensional cohomology group $H^1(E,\calO_E)$. The elliptic curve
$E$ is {\em supersingular\/} (or has {\em Hasse invariant\/ $0$}) if the map
$F$ above is zero, and is {\em ordinary\/} ({\em Hasse invariant\/ $1$})
otherwise. If $E=\Proj A$, then the map $F$ above is precisely the action of
the Frobenius on the socle of the injective hull of the residue field of $A$,
and hence $E$ is ordinary if and only if $A$ is an $F$-pure ring.

Let $f \in \ZZ[x_0,x_1,x_2]$ be a cubic polynomial defining a smooth elliptic
curve $E_\QQ \subset \PP^2_\QQ$. Then the Jacobian ideal of $f$ is
$(x_0,x_1,x_2)$-primary in $\QQ[x_0,x_1,x_2]$. Hence after localizing at an
appropriate nonzero integer $u$, the Jacobian ideal of $f$ in
$\ZZ[u^{-1}][x_0,x_1,x_2]$ contains high powers of $x_0,x_1$, and $x_2$.
Consequently, for all but finitely many prime integers $p$, the polynomial
$f \mod p$ defines a smooth elliptic curve $E_p \subset \PP^2_{\ZZ/p}$. If
the elliptic curve $E_\CC \subset \PP^2_\CC$ has complex multiplication,
then it is a classical result \cite{Deuring} that the {\em density\/} of the
supersingular prime integers $p$, i.e.,
$$
\lim_{n \to \infty}
\frac{|\{p \text{ prime}: p \le n \text{ and $E_p$ is supersingular}\}|}
{|\{p \text{ prime}: p \le n\}|}
$$
is $1/2$, and that this density is $0$ if $E_\CC$ does not have complex
multiplication. However, even if $E_\CC$ does not have complex multiplication,
the set of supersingular primes is infinite by \cite{Elkies}. It is conjectured
that if $E_\CC$ does not have complex multiplication, then the number of
supersingular primes less than $n$ grows asymptotically like
$C(\sqrt{n}/\log n)$, where $C$ is a positive constant, \cite{LT}.

Hartshorne and Speiser observed that the cohomological dimension of the
defining ideal of $E_p\times\PP^1_{\ZZ/p}$ varies with the prime $p$,
\cite[Example 3, p.\,75]{HS}. Their arguments use the notion of $F$-depth, and
we would like to point out how their results also follow from a recent theorem
of Lyubeznik:

\begin{theorem} \cite[Theorem 1.1]{Ly2} \label{dual}
Let $(R,\m)$ be a regular local ring containing a field of positive
characteristic, and $\a$ be an ideal of $R$. Then $H^i_\a(R)=0$ if and only if
there exists an integer $e \ge 1$ such that
$F^e:H^{\dim R -i}_\m(R/\a) \longrightarrow H^{\dim R -i}_\m(R/\a)$ is the zero
map, where $F^e$ denotes the $e$-th iteration of the Frobenius morphism.
\end{theorem}

\begin{cor}
Let $f \in \ZZ[x_0,x_1,x_2]$ be a cubic polynomial defining a smooth elliptic
curve $E_\QQ \subset \PP^2_\QQ$, and let
$\a \subset R=\ZZ[z_{ij}: 0\le i \le 2, \ 0 \le j \le n]$ be the ideal
defining the Segre embedding $E \times \PP^n \subset \PP^{3n+2}$. Then
$$
\cd(R/pR,\a)=\begin{cases}
2n+1 & \text{if} \quad \text{$E_p$ is supersingular,} \\
3n+1 & \text{if} \quad \text{$E_p$ is ordinary.}
\end{cases}
$$
\end{cor}

\begin{proof}
The ring $R/(\a+pR)$ may be identified with the Segre product $A\# B$ where
$$
A=\ZZ/p \ZZ[x_0,x_1,x_2]/(f) \quad \text{and} \quad B=\ZZ/p \ZZ [y_0,\dots,y_n].
$$
Let $p$ be a prime for which $E_p$ is smooth, in which case the ring
$A \otimes_{\ZZ/p} B$ and hence its direct summand $A\#B$ are normal. For
$k \ge 1$ the K\"unneth formula gives us
\begin{multline*}
H^{k+1}_\m(R/(\a+pR))=\bigoplus_{r \in \ZZ}
H^k(E_p \times \PP^n_{\ZZ/p},\calO_{E_p \times \PP^n_{\ZZ/p}}(r)) \\
=\bigoplus_{r \in \ZZ, i+j=k}
H^i(E_p,\calO_{E_p}(r)) \otimes H^j(\PP^n_{\ZZ/p},\calO_{\PP^n_{\ZZ/p}}(r)).
\end{multline*}
Hence
$$
H^{k+1}_\m(R/(\a+pR))=\begin{cases}
\ZZ/p \ZZ & \text{if} \quad k=1, \\
0 & \text{if} \quad 2 \le k \le n,
\end{cases}
$$
and the Frobenius action on the one-dimensional vector space $H^2_\m(R/(\a+pR))$
may be identified with the Frobenius
$$
H^1(E_p,\calO_{E_p}) \otimes H^0(\PP^n_{\ZZ/p},\calO_{\PP^n_{\ZZ/p}})
\overset{F}\longrightarrow
H^1(E_p,\calO_{E_p}) \otimes H^0(\PP^n_{\ZZ/p},\calO_{\PP^n_{\ZZ/p}}),
$$
which is the zero map precisely when $E_p$ is supersingular. Consequently every
element of $H^2_\m(R/(\a+pR))$ is killed by the Frobenius (equivalently, by an
iteration of the Frobenius) if and only if $E_p$ is supersingular. The
assertion now follows from Theorem~\ref{dual}.
\end{proof}

\begin{example}
The cubic polynomial $x^3+y^3+z^3$ defines a smooth elliptic curve $E_p$ in any
characteristic $p \neq 3$. It is easily seen that $E_p$ is supersingular for
primes $p \equiv 2 \mod 3$, and is ordinary if $p \equiv 1 \mod 3$. Let
$R=\ZZ[u,v,w,x,y,z]$. The defining ideal of $E \times \PP^1$ is the ideal $\a$
of $R$ generated by
\begin{multline*}
u^3+v^3+w^3, \ u^2x+v^2y+w^2z, \ ux^2+vy^2+wz^2, \ x^3+y^3+z^3, \\
vz-wy, \ wx-uz, \ uy-vx.
\end{multline*}
If $p \neq 3$ is a prime integer, then $H^4_\a(R/pR)=0$ if and only if 
$p \equiv 2 \mod 3$, and consequently
$$
\cd(R/pR,\a)=\begin{cases}
3 & \text{if} \quad p \equiv 2 \mod 3, \\
4 & \text{if} \quad p \equiv 1 \mod 3.
\end{cases}
$$
\end{example}

As the above example shows, the cohomological dimension $\cd(R/pR,\a)$ varies
with the characteristic $p$, so we cannot use local cohomology to complete the
proof of Theorem~\ref{thm-ara} in arbitrary prime characteristic. We instead
use \'etale cohomology to show that the defining ideal of {\em any\/}
projective embedding of $E \times \PP^1$ cannot be a set-theoretic complete
intersection which, in particular, completes the proof of
Theorem~\ref{thm-ara}.

\begin{theorem}\label{etale}
Let $E$ be a smooth elliptic curve. Then the defining ideal of any projective
embedding of $E \times \PP^1$ is not a set-theoretic complete intersection.
\end{theorem}

\begin{proof}
Consider an embedding $E\times \PP^1 \subset \PP^k$. Then the defining ideal
$\a$ has height $k-2$. We need to prove that $\ara \a > k-2$, and for this we
may replace the field $K$ by its separable closure. Let $\ell$ be a prime
integer different from the characteristic of $K$. We shall use the \'etale
cohomology groups $H^i_\et(-;\ZZ/\ell\ZZ)$, i.e., with coefficients in
$\ZZ/\ell\ZZ$.

If $\a=\rad(g_1,\dots,g_{k-2})$, then the complement $U$ of $X=E\times \PP^1$
in $\PP^k$ can be covered by the affine open sets $U_j=D_+(g_j)\subset\PP^k$
for $1 \le j \le k-2$. Each $U_j$ is an affine smooth variety of dimension $k$,
and so $H^i_\et(U_j;\ZZ/\ell\ZZ)=0$ for all $i>k$ by \cite[Thm.~VI.7.2]{Milne}.
The Mayer-Vietoris principle \cite[III.2.24]{Milne} now implies that
$H^{2k-2}_\et(U;\ZZ/\ell\ZZ)=0$. We shall show that this leads to a
contradiction.

Since $H^1_\et(E;\ZZ/\ell\ZZ)$ is nonzero, the K\"unneth formula
\cite[VI.8.13]{Milne} implies that $H^1_\et(X;\ZZ/\ell\ZZ)$ is nonzero as well.
Since $X$ is proper, the first compactly supported \'etale cohomology of $X$ is
$H^1_{\et,c}(X;\ZZ/\ell\ZZ)=H^1_\et(X;\ZZ/\ell\ZZ)$, \cite[III.1.29]{Milne}.
By \cite[(1.4a)]{Ly-etale} and \cite[Cor.~VI.11.2]{Milne} there is a natural
isomorphism between $H^i_{\et,c}(X;\ZZ/\ell\ZZ)$ and the dual of
$H^{2k-i}_{\et,X}(\PP^k;\ZZ/\ell\ZZ)$, so it follows that
$H^{2k-1}_{\et,X}(\PP^k;\ZZ/\ell\ZZ)$ is nonzero. By \cite[III.1.25]{Milne} we
have an exact sequence
$$
\begin{CD}
H^i_\et(U;\ZZ/\ell\ZZ) @>>> H^{i+1}_{\et,X}(\PP^k;\ZZ/\ell\ZZ)
@>>> H^{i+1}_\et(\PP^k;\ZZ/\ell\ZZ).
\end{CD}
$$
But $H^{2k-1}_\et(\PP^5;\ZZ/\ell\ZZ)=0$ by \cite[VI.5.6]{Milne}, which gives a
contradiction.
\end{proof}

\section{The set-theoretically Cohen-Macaulay property}

Given an affine variety $V$, it is an interesting question whether $V$ supports
a Cohen-Macaulay scheme, i.e., whether there exists a Cohen-Macaulay ring $R$
such that $V$ is isomorphic to $\Spec R_\red$. More generally, let $R$ be a
regular local ring. We say that an ideal $\a\subset R$ is {\em
set-theoretically Cohen-Macaulay\/} if there exists an ideal $\b \subset R$
with $\rad\b=\rad\a$ for which the ring $R/\b$ is Cohen-Macaulay. A homogeneous
ideal $\a$ of a polynomial ring $R$ is {\em set-theoretically Cohen-Macaulay\/}
if $\a R_\m$ is a set-theoretically Cohen-Macaulay ideal of $R_\m$, where $\m$
is the homogeneous maximal ideal of $R$.

There is a well-known example of a determinantal ideal which is not a
set-theoretic complete intersection, but is Cohen-Macaulay (and hence
set-theoretically Cohen-Macaulay). For an integer $n\ge 2$, let $X=(x_{ij})$ be
an $n \times (n+1)$ matrix of variables over a field $K$, and let $R$ be the
localization of the polynomial ring $K[x_{ij}:1\le i\le n,1\le j\le n+1]$ at
its homogeneous maximal ideal. Let $\a$ be the ideal of $R$ generated by the $n
\times n$ minors of the matrix $X$. If $K$ has characteristic zero, then
$H^{n+1}_\a(R)$ is nonzero by an argument due to Hochster, \cite[Remark
3.13]{HL}, so $\ara \a=n+1$. If $K$ has positive characteristic it turns out
that $H^{n+1}_\a(R)=0$, but nevertheless the ideal $\a$ has arithmetic rank
$n+1$, \cite{Newstead, BS}. In particular, $\a$ is not a set-theoretic complete
intersection though $R/\a$ is Cohen-Macaulay.

We next show that for a smooth elliptic curve $E_\QQ \subset \PP^2_\QQ$, the
defining ideal of $E_\QQ \times \PP^1_\QQ \subset \PP^5_\QQ$ is not
set-theoretically Cohen-Macaulay. We begin with a lemma of Huneke, \cite[page
599]{EMS}. We include a proof here for the convenience of the reader.

\begin{lemma}[Huneke]\label{setcm}
Let $\a$ be an ideal of a regular local ring $R$ of characteristic $p>0$. If
the ring $R/\a$ is $F$-pure and not Cohen-Macaulay, then the ideal $\a$ is not
set-theoretically Cohen-Macaulay.
\end{lemma}

\begin{proof}
Note that $\a$ is a radical ideal since $R/\a$ is $F$-pure. Let $\b$ be an ideal
of $R$ with $\rad \b=\a$, and choose $x_1,\dots, x_d \in R$ such that their
images form a system of parameters for $R/\a$ and $R/\b$. Since $R/\a$ is is
not Cohen-Macaulay, there exist $k\in\NN$ and $y \in R$ such that
$$
yx_k \in (x_1,\dots,x_{k-1})R+\a \quad \text{and} \quad
y \notin (x_1,\dots,x_{k-1})R+\a.
$$
Let $q=p^e$ be a prime power such that $\a^{[q]} \subseteq \b$. Then
$$
y^qx_k^q \in (x_1^q,\dots,x_{k-1}^q)R+\a^{[q]} \ \subseteq \
(x_1^q,\dots,x_{k-1}^q)R+\b.
$$
If $R/\b$ is Cohen-Macaulay, then
$$
y^q \in (x_1^q,\dots,x_{k-1}^q)R+\b \ \subseteq \
(x_1^q,\dots,x_{k-1}^q)R+\a.
$$
The hypothesis that $R/\a$ is $F$-pure implies that
$y \in (x_1,\dots,x_{k-1})R+\a$, which is a contradiction.
\end{proof}

For the remainder of this section, $R_\ZZ$ will denote a polynomial ring over
the integers, and we use the notation $R_\QQ=R_\ZZ\otimes_\ZZ\QQ$ and
$R_p=R_\ZZ\otimes_\ZZ\ZZ/p\ZZ$.

\begin{lemma}\label{modp}
Let $\a$ be an ideal of $R_\ZZ=\ZZ[z_1,\dots,z_m]$, and consider the
multiplicative set $W=R_\ZZ \setminus (z_1,\dots,z_m)R_\ZZ$. If
$W^{-1}R_\QQ/ \a W^{-1} R_\QQ$ is Cohen-Macaulay, then the rings
$W^{-1}R_p/ \a W^{-1} R_p$ are Cohen-Macaulay for all but finitely many prime
integers $p$.
\end{lemma}

\begin{proof}
Let $x_1,\dots,x_d \in R_\ZZ$ be elements whose images form a system of
parameters for $W^{-1}R_\QQ/\a W^{-1}R_\QQ$. Since this ring is Cohen-Macaulay,
there exists an element $f$ in the multiplicative set $W$ such that
$x_1,\dots,x_d$ is a regular sequence on $S/\a S$ where $S=R_\ZZ[f^{-1}]$.
Moreover, we may choose $f$ in such a way that
$$
z_i \in \rad(\a S+(x_1,\dots,x_d)S), \qquad 1\le i\le m.
$$
These conditions are preserved if we enlarge the ring $S$ by inverting finitely
many nonzero integers. By the result on generic freeness,
\cite[Theorem 24.1]{Matsumura}, we may assume (after replacing $f$ by a nonzero
integer multiple $uf$ and $S$ by its localization at the element $u$) that 
$S$, $S/\a$, and each of
$$
\frac{S}{\a S+(x_1,\dots,x_i)S}, \qquad 1 \le i \le d,
$$
are free $\ZZ[u^{-1}]$-modules. In particular, for all $1\le i\le d-1$, we have
short exact sequences of free $\ZZ[u^{-1}]$-modules,
$$
\begin{CD}
0 @>>> \frac{S}{\a S+(x_1,\dots,x_i)S} @>{x_{i+1}}>>
\frac{S}{\a S+(x_1,\dots,x_i)S} @>>> \frac{S}{\a S+(x_1,\dots,x_{i+1})S}
@>>> 0.
\end{CD}
$$
Let $p$ be any prime integer not dividing $u$, and apply
$(-)\otimes_{\ZZ[u^{-1}]} \ZZ/p \ZZ$ to the sequences above. The resulting exact
sequences show that $x_1,\dots,x_d$ is a regular sequence on 
$$
\frac{S}{\a S}\otimes_{\ZZ[u^{-1}]} \frac{\ZZ}{p \ZZ} \cong
\frac{R_p[f^{-1}]}{\a R_p[f^{-1}]},
$$
and hence on $W^{-1}R_p/ \a W^{-1} R_p$ as required.
\end{proof}

\begin{theorem}
Let $E_\QQ \subset \PP^2_\QQ$ be a smooth elliptic curve. Then the defining
ideal of the Segre embedding $E_\QQ \times \PP^1_\QQ \subset \PP^5_\QQ$ is not
set-theoretically Cohen-Macaulay.
\end{theorem}

\begin{proof}
Let $\a\subset R_\ZZ=\ZZ[z_{ij}: 0\le i\le 2,\ 0\le j\le 1]$ be an ideal
such that $\a R_\QQ\subset R_\QQ$ is the defining ideal of
$E_\QQ\times\PP^1_\QQ\subset\PP^5_\QQ$. There exist infinitely many prime
integers $p$ such that $E_p\subset\PP^2_{\ZZ/p}$ is a smooth ordinary
elliptic curve. For these infinitely many primes, the ring $R_p/\a R_p$ is
$F$-pure and not Cohen-Macaulay, and hence the ideal $\a R_p$ is not
set-theoretically Cohen-Macaulay by Lemma~\ref{setcm}. It now follows from
Lemma~\ref{modp} that the ideal $\a R_\QQ \subset R_\QQ$ is not
set-theoretically Cohen-Macaulay.
\end{proof}

\begin{remark}
If $E_p\subset\PP^2_{\ZZ/p}$ is an ordinary elliptic curve, then the defining
ideal of the Segre embedding $E_p\times\PP^1_{\ZZ/p}\subset\PP^5_{\ZZ/p}$ is
not set-theoretically Cohen-Macaulay by Lemma~\ref{setcm}, since the
corresponding homogeneous coordinate ring is $F$-pure and not Cohen-Macaulay.
If $E$ is supersingular, we do not know whether the defining ideal is
set-theoretically Cohen-Macaulay.
\end{remark}

\begin{remark}
Lemma~\ref{setcm} can be strengthened by Lyubeznik's theorem as follows: Let
$(R,\m)$ be a regular local ring of positive characteristic, and $\a$ and $\b$
be ideals of $R$ such that $R/\a$ is $F$-pure and $\rad \b=\a$. If the local
cohomology module $H^i_\m(R/\a)$ is nonzero for some integer $i$, then
Theorem~\ref{dual} implies that
$$
H^{\dim R - i}_\a(R)=H^{\dim R - i}_\b(R)
$$
is nonzero as well. Using Theorem~\ref{dual} once again, it follows that
$H^i_\m(R/\b)$ is nonzero. This implies in particular that
$$
\depth R/\a \ge \depth R/\b.
$$
Consequently if $\b$ is an ideal of a regular local ring $R$ of positive
characteristic such that $R/\rad \b$ is $F$-pure, then
$$
\depth R/\rad \b \ge \depth R/\b.
$$
This is false without the assumption that $R/\rad \b$ is $F$-pure: let
$R=K[[w,x,y,z]]$ where $K$ is a field of positive characteristic, and consider
the ideal
$$
\p=(xy-wz, y^3-xz^2, wy^2-x^2z, x^3-w^2y) \subset R.
$$
Then $R/\p \cong K[[s^4, s^3t, st^3, t^4]]$ which is not Cohen-Macaulay.
Hartshorne proved that $\p$ is a set-theoretic complete intersection, i.e.,
that $\p=\rad(f,g)$ for elements $f,g \in R$, \cite{Hartshorne}. Hence, in
this example,
$$
\depth R/\rad(f,g)=1 < \depth R/(f,g)=2.
$$
\end{remark}

\section*{Acknowledgments}
It is a pleasure to thank Gennady Lyubeznik for several helpful discussions.


\begin{thebibliography}{10}

\bibitem{Bott-Tu} R.~Bott and L.~W.~Tu, \emph{Differential forms in algebraic
topology}, Graduate Texts in Mathematics, vol.~82, Springer-Verlag, New York,
1982.

\bibitem{BS} W.~Bruns and R.~Schw{\"a}nzl, \emph{The number of equations
defining a determinantal variety}, Bull. London Math. Soc. \textbf{22} (1990),
no.~5, 439--445.

\bibitem{Deuring} M.~Deuring, \emph{Die {T}ypen der {M}ultiplikatorenringe
elliptischer {F}unktionenk\"orper}, Abh. Math. Sem. Hansischen Univ.
\textbf{14} (1941), 197--272.

\bibitem{EMS} D.~Eisenbud, M.~Musta{\c{t}}{\v{a}}, and M.~Stillman,
\emph{Cohomology on toric varieties and local cohomology with monomial
supports}, J. Symbolic Comput. \textbf{29} (2000), no.~4-5, 583--600, Symbolic
computation in algebra, analysis, and geometry (Berkeley, CA, 1998).

\bibitem{Elkies} N.~D.~Elkies, \emph{The existence of infinitely many
supersingular primes for every elliptic curve over {${\bf Q}$}}, Invent. Math.
\textbf{89} (1987), no.~3, 561--567.

\bibitem{Hartshorne} R.~Hartshorne, \emph{Complete intersections in
characteristic {$p>0$}}, Amer. J. Math. \textbf{101} (1979), no.~2, 380--383.

\bibitem{HS} R.~Hartshorne and R.~Speiser, \emph{Local cohomological dimension
in characteristic {$p$}}, Ann. of Math. (2) \textbf{105} (1977), no.~1, 45--79.

\bibitem{HR} M.~Hochster and J.~L.~Roberts, \emph{The purity of the {F}robenius
and local cohomology}, Advances in Math. \textbf{21} (1976), no.~2, 117--172.

\bibitem{HL} C.~Huneke and G.~Lyubeznik, \emph{On the vanishing of local
cohomology modules}, Invent. Math. \textbf{102} (1990), no.~1, 73--93.

\bibitem{Iversen} B.~Iversen, \emph{Cohomology of sheaves}, Universitext,
Springer-Verlag, Berlin, 1986.

\bibitem{LT} S.~Lang and H.~Trotter, \emph{Frobenius distributions in {${\rm
GL}\sb{2}$}-extensions}, Springer-Verlag, Berlin, 1976, Distribution of
Frobenius automorphisms in ${\rm GL}\sb{2}$-extensions of the rational numbers,
Lecture Notes in Mathematics, Vol. 504.

\bibitem{Ly1} G.~Lyubeznik, \emph{The number of defining equations of affine
algebraic sets}, Amer. J. Math. \textbf{114} (1992), no.~2, 413--463.

\bibitem{Ly-etale} \bysame, \emph{\'{E}tale cohomological dimension and the
topology of algebraic varieties}, Ann. of Math. (2) \textbf{137} (1993), no.~1,
71--128.

\bibitem{Ly2} \bysame, \emph{On the vanishing of local cohomology in
characteristic $p>0$}, Preprint (2004).

\bibitem{Matsumura} H.~Matsumura, \emph{Commutative ring theory}, second ed.,
Cambridge Studies in Advanced Mathematics, vol.~8, Cambridge University Press,
Cambridge, 1989.

\bibitem{Milne} J.~S. Milne, \emph{\'{E}tale cohomology}, Princeton
Mathematical Series, vol.~33, Princeton University Press, Princeton, N.J.,
1980.

\bibitem{Newstead} P.~E.~Newstead, \emph{Some subvarieties of {G}rassmannians
of codimension {$3$}}, Bull. London Math. Soc. \textbf{12} (1980), no.~3,
176--182.
\end{thebibliography}
\end{document}